\documentclass[12pt,twoside]{article}
\usepackage[english]{babel}
\usepackage{umj_3_2}           % Style package of Ural Mathematical Journal
\usepackage{graphicx, srcltx}   % For graphics
\usepackage{cite}               % For correct bibliography
\usepackage{url}                                % use command \cite{1,2,3,4,8,9,10}
\usepackage{amsmath}            % If the package features are used
\usepackage{amsfonts,amssymb}   % If there is need to use math alphabets
\usepackage{srcltx}             % For easy move between WiNEdT and YAP

\newcommand{\Leftclt}{ Victor Nijimbere}   % Left page header - name & family name of article author
\newcommand{\Rightclt}{The non-elementary integral $\int e^{\lambda x^\alpha} dx$, $\alpha\ge2 $, and other related integrals} %% Right page header - Short title of article
\usepackage{authblk}

\begin{document}

\date{}

\Mainclt            % Page header print on the first page

% Enter article title and authors
% Reference to grants is a footnote to article title. You may use also Acknowledgments section for this purpose if you wish.

\begin{Titul}
{\large \bf  EVALUATION OF THE NON-ELEMENTARY INTEGRAL\\[1ex]  $\int e^{\lambda x^\alpha} dx$, $\ \alpha\ge2$, AND OTHER RELATED INTEGRALS}\\ \vspace{5mm}
%\footnote{This work was supported by ... (project no.NN-N).
%You may use also Acknowledgments section for this purpose if you wish.}\\[3ex]
{{\bf Victor Nijimbere} \\ [2ex] {\small School of Mathematics and Statistics, Carleton University,\\ Ottawa, Ontario, Canada\\ victornijimbere@gmail.com}\\[3ex]}
%{\bf  Second B. Author}\\[2ex] {\small University, City, Country, Email }\\[5ex]}
\end{Titul}

\begin{Anot}
{\bf Abstract:} A formula for the non-elementary integral $\int e^{\lambda x^\alpha} dx$ where $\alpha$ is real and greater or equal two,
is obtained in terms of the confluent hypergeometric function $_{1}F_1$ by expanding the integrand as a Taylor series.
This result is verified by directly evaluating the area under the Gaussian Bell curve, corresponding to $\alpha=2$,
using the asymptotic expression for the confluent hypergeometric function and the Fundamental Theorem of Calculus (FTC).
Two different but equivalent expressions, one in terms of the confluent hypergeometric
function $_{1}F_1$ and another one in terms of the hypergeometric function $_1F_2$, are obtained for each of these integrals,
$\int\cosh(\lambda x^\alpha)dx$, $\int\sinh(\lambda x^\alpha)dx$, $\int \cos(\lambda x^\alpha)dx$ and $\int\sin(\lambda x^\alpha)dx$,
$\lambda\in \mathbb{C},\alpha\ge2$. And the hypergeometric function $_1F_2$ is expressed in terms of the confluent hypergeometric function $_1F_1$. Some of the applications of the non-elementary integral $\int e^{\lambda x^\alpha} dx, \alpha\ge 2$ such as the Gaussian distribution and the Maxwell-Bortsman distribution are given.

\end{Anot}

%%%%%%%%%%%%%%%%%%%%%%%%%%%%%%%%%%%%%%%%%%%%
 %Please replace what follows by the body %%
 %of your article                         %%
                                          %%
%%%%%%%%%%%%%%%%%%%%%%%%%%%%%%%%%%%%%%%%%%%%

\section{Introduction}\label{sec:1}
\setcounter{equation}{0}

%This article illustrates the preparation of a paper for Ural Mathematical Journal using \LaTeX2e.

%Your paper must be prepared in {\bf English}. To ensure that your paper will be reproduced
%clearly and in the proper size and form, please observe the following instructions.

%The length of the title of your paper should not exceed two lines.
%The paper should begin with an abstract, followed by key words, and then an introduction.
%At the end, it should feature a conclusion and a list of references.
%The maximum length of your paper, including figures and tables,should be 10-20 pages.

%In addition to length and
%format restrictions given here, all papers will be reviewed and must meet technical standards.

%\begin{verbatim}\begin{assen}...\end{assen}\end{verbatim}
%For example,
%\begin{teoen}[{\cite[Theorem~2]{DeVore-Lorentz, Geronimus-1936}}]
%Theorem text goes here.
%\end{teoen}
%\begin{teoen}[(Eremenko, Yuditskii {\rm \cite[p.~25]{Eremenko-Yuditskii}})]
%Theorem text goes here.
%\end{teoen}
%\proofen    Proof text goes here.  \hfill$\square$\\[1ex]%--- P r o o f.
%\proofnpen follows from $\ldots$\,... \hfill$\square$\\[1ex]% --- P r o o f

%\emph{Remark~1.} Text...\\[1ex] % Remark (Çàìå÷àíèå)
%\emph{Example~1.} Text...\\[1ex] % Example

%We use \hfill$\square$ for the end of the proof.

\begin{defen} An elementary function is a function of one variable built up
using that variable and constants, together with a finite number of repeated
algebraic operations and the taking of exponentials and logarithms\emph{ \cite{6}}.
\label{defen:1}
\end{defen}

In 1835, Joseph Liouville established conditions in his theorem, known as Liouville 1835's Theorem \cite{4,6}, which can be used to determine whether an indefinite integral is elementary or non-elementary. Using Liouville 1835's Theorem, one can show that the indefinite integral $\int e^{\lambda x^\alpha}dx,$ ${\alpha \ge 2}$, is non-elementary \cite{4}, and to my knowledge, no one has evaluated this non-elementary integral before.

For instance, if $\alpha = 2$, $\lambda = -\beta^2 < 0$, where $\beta$ is a real constant, the area
under the Gaussian Bell curve can be calculated using double integration and then polar coordinates to obtain

\begin{equation}
\int\limits_{-\infty}^{+\infty} e^{-\beta^2 x^2}dx =\frac{\sqrt{\pi}}{\beta}.
\label{eq1}
\end{equation}
Is that possible to evaluate (\ref{eq1}) by directly using the Fundamental Theorem of
Calculus (FTC) as in equation (\ref{eq2})?
\begin{equation}
\int\limits_{-\infty}^{+\infty} e^{-\beta^2 x^2}dx=
\lim_{t\to -\infty}\int\limits_{t}^{0} e^{-\beta^2 x^2}dx+\lim_{t\to +\infty}\int\limits_{0}^{t} e^{-\beta^2 x^2}dx.
\label{eq2}
\end{equation}

The Central limit Theorem (CLT) in Probability theory \cite{2} states that the
probability that a random variable $x$ does not exceed some observed value $z$
\begin{equation}
P(X<z)=\frac{1}{\sqrt{2\pi}}\int\limits_{-\infty}^{z} e^{-\frac{x^2}{2}}dx.
\label{eq3}
\end{equation}
So if we know the antiderivative of the function $g(x) = e^{\lambda x^2}$, we may choose
to use the FTC to calculate the cumulative probability $P(X < z)$ in (\ref{eq3}) when the
value of $z$ is given or is known, rather than using numerical integration.

The Maxwell-Boltsman distribution in gas dynamics,
\begin{equation}
F(v)=\theta\int\limits_{0}^{v} x^2e^{-\gamma x^2}dx,
\label{eq4}
\end{equation}
where $\theta$ and $\gamma$ are some positive constants that depend on the properties of the gas and $v$ is the gas speed, is another application.

There are many other examples where the antiderivative of $g(x) = e^{\lambda x^\alpha}$, $\alpha \ge
2$ can be useful. For example, using the FTC, formulas for integrals such as

\begin{equation}
\int\limits_{x}^{\infty} e^{t^{2n+1}}dt, x<\infty;\quad \int\limits_{x}^{\infty}e^{-t^{2n+1}}dt, x > -\infty; \quad
\int\limits_{x}^{\infty}t^{2n}e^{-t^2}dt, x\le\infty,
\label{eq5}
\end{equation}
where $n$ is a positive integer, can be obtained if the antiderivative of $g(x)=e^{\lambda x^\alpha}$, $\alpha \ge 2$ is known.

In this paper, the antiderivative of $g(x)=e^{\lambda x^\alpha},$ $\alpha \ge 2$, is expressed in terms of a special function, the confluent hypergeometric $_1F_1$ \cite{1}. And the confluent hypergeometric $_1F_1$ is an entire function \cite{3}, and its properties are well known \cite{1,5}.
The main goal here is to consider the most general case with $\lambda$ complex $(\lambda \in \mathbb{C})$,
evaluate the non-elementary integral $\int e^{\lambda x^\alpha}$, $\alpha\ge2$ and thus make
possible the use of the FTC to compute the definite integral

\begin{equation}
\int\limits_{A}^{B} e^{\lambda x^\alpha}dx,
\label{eq6}
\end{equation}
for any A and B. And once (\ref{eq6}) is evaluated, then integrals such as (\ref{eq1}), (\ref{eq2}),
(\ref{eq3}), (\ref{eq4}) and (\ref{eq5}) can also be evaluated using the FTC.

Using the hyperbolic and Euler identities, \begin{gather*}\cosh(\lambda x^\alpha) = (e^{\lambda x^\alpha}+ e^{-\lambda x^\alpha})/2, \quad
\sinh(\lambda x^\alpha) = (e^{\lambda x^\alpha}-e^{-\lambda x^\alpha})/2,\\
\cos(\lambda x^\alpha) = (e^{i\lambda x^\alpha}+ e^{-i\lambda x^\alpha})/2, \quad
\sin(\lambda x^\alpha) = (e^{i\lambda x^\alpha}- e^{-i\lambda x^\alpha})/(2i),
\end{gather*}
 the integrals
%$\int \cosh(\lambda x^\alpha) dx, \int \sinh(\lambda x^\alpha) dx, \int \cos(\lambda x^\alpha) dx$ and $\int \sin(\lambda x^\alpha) dx$,
\begin{equation}
\int \cosh(\lambda x^\alpha) dx,\quad \int \sinh(\lambda x^\alpha) dx,\quad \int \cos(\lambda x^\alpha) dx \quad \hspace{.12cm}\mbox{and}\hspace{.12cm} \quad
 \int \sin(\lambda x^\alpha) dx, \alpha \ge 2,
\label{eq7}
\end{equation}
are evaluated in terms of $_1F_1$ for any constant $\lambda$. They
are also expressed in terms of the hypergeometric $_1F_2$. And some expressions
of the hypergeometric function $_1F_2$ in terms of the confluent hypergeometric
function $_1F_1$ are therefore obtained.

For reference, we shall first define the confluent confluent hypergeometric function $_1 F_1$ and the hypergeometric  function $_1 F_2$  before we proceed to the main aims of this paper (see sections~\ref{sec:2} and~\ref{sec:3}).

\begin{defen}The confluent hypergeometric function, denoted as $_1F_1$, is a special function given by the series \emph{\cite{1,5}}
\begin{equation}
_1 F_1(a;b; x)=\sum\limits_{n=0}^{\infty}\frac{(a)_n}{(b)_n}\frac{x^n}{n!},
\label{confluent}
\end{equation}
where $a$ and $b$ are arbitrary constants, $(\vartheta)_n=\Gamma(\vartheta+n)/\Gamma(\vartheta)$ \emph{(}Pochhammer's notation \emph{\cite{1}}\emph{)} for any complex $\vartheta$, with $(\vartheta)_0=1$, and $\Gamma$ is the standard gamma function \emph{\cite{1}}.
\label{defen:2}
\end{defen}

\begin{defen} The hypergeometric function $_1F_2$ is a special function given by the series \emph{\cite{1,5}}
\begin{equation}
_1 F_2(a;b,c; x)=\sum\limits_{n=0}^{\infty}\frac{(a)_n}{(b)_n (c)_n}\frac{x^n}{n!},
\label{confluent}
\end{equation}
where $a,b$ and $c$ are arbitrary constants, and $(\vartheta)_n=\Gamma(\vartheta+n)/\Gamma(\vartheta)$ \emph{(}Pochhammer's notation \emph{\cite{1}}\emph{)} as in Definition~\ref{defen:2}.
\label{defen:3}
\end{defen}

\section{Evaluation of $\int_A^B e^{\lambda x^\alpha}dx$}\label{sec:2}
\setcounter{equation}{0}
%Symbols and acronyms should be defined the first time they appear. Use the International System (SI) of units.

\begin{prpen}The function $G(x) = x\ _1F_1\left(\frac{1}{\alpha}; \frac{1}{\alpha}+1; \lambda x^\alpha\right)$, where $_1F_1$ is a confluent hypergeometric function \emph{\cite{1}}, $\lambda$ is an arbitrarily constant and $\alpha\ge 2$,
is the antiderivative of the function $g(x) = e^{\lambda x^\alpha}$. Thus,
\begin{equation}
\int e^{\lambda x^\alpha}dx=x\ _1F_1\left(\frac{1}{\alpha}; \frac{1}{\alpha}+1; \lambda x^\alpha\right)+C.
\label{eq8}
\end{equation}
\label{prpen:1}
\end{prpen}
\proofen   We expand $g(x) = e^{\lambda x^\alpha}$ as a Taylor series and integrate the series
term by term. We also use the Pochhammer's notation \cite{1} for the gamma
function, $\Gamma(a + n) = \Gamma(a)(a)_n$, where $(a)_n = a(a + 1) \cdots (a + n - 1)$,
and the property of the gamma function $\Gamma(a + 1) = a\Gamma(a)$ \cite{1}. For example,
$\Gamma(n + a + 1) = (n + a)\Gamma (n + a)$. We then obtain
\begin{align}\begin{aligned}
\int g(x) dx &=\int e^{\lambda x^\alpha}dx=\sum\limits_{n=0}^{\infty}\frac{\lambda^n}{n!}\int x^{\alpha n}dx  \\&
=\sum\limits_{n=0}^{\infty}\frac{\lambda^n}{n!}\frac{x^{\alpha n+1}}{\alpha n+1}+C
=\frac{x}{\alpha}\sum\limits_{n=0}^{\infty}\frac{(\lambda x^{\alpha})^n}{\left(n+\frac{1}{\alpha}\right)n!}+C \\&=\frac{x}{\alpha}\sum\limits_{n=0}^{\infty}\frac{\Gamma\left(n+\frac{1}{\alpha}\right)}{\Gamma\left(n+\frac{1}{\alpha}+1\right)}\frac{(\lambda x^{\alpha})^n}{n!}+C
 \\&
={x}\sum\limits_{n=0}^{\infty}\frac{\left(\frac{1}{\alpha}\right)_n}{\left(\frac{1}{\alpha}+1\right)_n}\frac{(\lambda x^{\alpha})^n}{n!}+C
 \\&= x\ _1F_1\left(\frac{1}{\alpha}; \frac{1}{\alpha}+1; \lambda x^\alpha\right)+C=G(x) + C. \quad \square
\label{eq9}
\end{aligned}\end{align}
%\hfill$\square$\\[1ex]
%\begin{gather}\begin{gathered}
%\int g(x) dx =\int e^{\lambda x^\alpha}dx=\sum\limits_{n=0}^{\infty}\frac{\lambda^n}{n!}\int x^{\alpha n}dx  \\
%=\sum\limits_{n=0}^{\infty}\frac{\lambda^n}{n!}\frac{x^{\alpha n+1}}{\alpha n+1}+C
%=\frac{x}{\alpha}\sum\limits_{n=0}^{\infty}\frac{(\lambda x^{\alpha})^n}{\left(n+\frac{1}{\alpha}\right)n!}+C \\=\frac{x}{\alpha}\sum\limits_{n=0}^{\infty}\frac{\Gamma\left(n+\frac{1}{\alpha}\right)}{\Gamma\left(n+\frac{1}{\alpha}+1\right)}\frac{(\lambda x^{\alpha})^n}{n!}+C
%\\
%={x}\sum\limits_{n=0}^{\infty}\frac{\left(\frac{1}{\alpha}\right)_n}{\left(\frac{1}{\alpha}+1\right)_n}\frac{(\lambda x^{\alpha})^n}{n!}+C
% \\ = x\ _1F_1\left(\frac{1}{\alpha}; \frac{1}{\alpha}+1; \lambda x^\alpha\right)+C=G(x) + C. \quad\square
%\label{eq9}
%\end{gathered}\end{gather}

%\\[1ex]%--- P r o o f.

%\begin{exaen}

\medskip
\emph{Example~1.}
We can now evaluate $\int x^{2n}e^{\lambda x^2}dx$ in terms of the confluent hypergeometric function. Using integration by parts,
\begin{equation}
\int x^{2n}e^{\lambda x^2}dx=\frac{ x^{2n-1}}{2\lambda}e^{\lambda x^2}-\frac{2n-1}{2\lambda}\int x^{2n-2}e^{\lambda x^2}dx.
\label{eq10}
\end{equation}
\begin{enumerate}
\item For instance, for $n=1$,
\begin{equation}
\int x^{2}e^{\lambda x^2}dx=\frac{ x}{2\lambda}e^{\lambda x^2}-\frac{1}{2\lambda}\int e^{\lambda x^2}dx=\frac{ x}{2\lambda}e^{\lambda x^2}-\frac{x}{2\lambda}\ _1F_1\left(\frac{1}{2}; \frac{3}{2}; \lambda x^2\right)+C.
\label{eq10}
\end{equation}
\item For $n=2$,
\begin{equation}
\int x^{4}e^{\lambda x^2}dx=\frac{ x^3}{2\lambda}e^{\lambda x^2}-\frac{3}{2\lambda}\int x^2e^{\lambda x^2}dx=\frac{ x^3}{2\lambda}e^{\lambda x^2}-\frac{ 3x}{4\lambda^2}e^{\lambda x^2}+\frac{3x}{4\lambda^2} \ _1F_1\left(\frac{1}{2}; \frac{3}{2}; \lambda x^2\right)+C.
\label{eq11}
\end{equation}
\end{enumerate}
\label{exaen:1}
%\end{exaen}

\medskip

%\begin{exaen}
\emph{Example~2.}
 Using the method of integrating factor, the first-order ordinary differential equation
\begin{equation}
y^\prime+2xy=1
\label{eq12}
\end{equation}
has solution
\begin{equation}
y(x)=e^{-x^2}\left(\int e^{x^2}dx+C\right)=xe^{-x^2}\ _1F_1\left(\frac{1}{2}; \frac{3}{2};  x^2\right)+Ce^{-x^2}.
\label{eq13}
\end{equation}
%\end{exaen}

\medskip

 Assuming that the function $G(x)$ (see Proposition \ref{prpen:1}) is unknown, in the following lemma, we use the properties of function $g(x)$ to establish the properties of $G(x)$ such as the inflection points and the behavior as $x\to\pm\infty$.

\begin{lemen}Let the function $G(x)$ be an antiderivative of $g(x)= e^{\lambda x^\alpha}, \lambda\in\mathbb{C}$ with $\alpha\ge2$.
\begin{enumerate}
\item If the real part of $\lambda$ is negative $(<0)$ and $\alpha$ is even, then the limits $lim_{x\to-\infty}G(x)$ and
$lim_{x\to+\infty}G(x)$ are finite (constants). And thus the Lebesgue integral $\int_{-\infty}^{\infty} |e^{\lambda x^\alpha}|dx<\infty$.

\item  If $\lambda$ is real $(\lambda\in\mathbb{R})$, then the point $(0, G(0)) = (0, 0)$ is an inflection point of the
curve $Y = G(x),$ $x\in\mathbb{R}$.

\item And if $\lambda\in\mathbb{R}$ and $\lambda<0$, and $\alpha$ is even, then the limits $\lim_{x\to-\infty}G(x)$ and
$\lim_{x\to+\infty}G(x)$ are finite. And there exists real constant $\theta>0$ such that
 limits $\lim_{x\to-\infty}G(x)=-\theta$ and $\lim_{x\to+\infty}G(x)=\theta$.
\end{enumerate}
\label{lemen:1}
\end{lemen}
\proofen
\begin{enumerate}
\item For complex $\lambda=\lambda_r+i\lambda_i$, where the subscript $r$ and $i$ stand for real and imaginary parts respectively, the function $g(x)=g(z)=e^{z^\alpha}$ where $z=(\lambda_r+i\lambda_i)^{1/\alpha}x,$ $ \alpha\ge2$, is an entire function on $\mathbb{C}$. And if $ \lambda_r<0$ and $\alpha$ is even implies $\mbox{Re}(z^\alpha)$ is always negative regardless of the values of $x$. And so, if $|z|\to\infty$ (or $x\to\pm\infty$), then $g(z)=0$ ($g(z)\to 0$) (or $g(x)=0$ as $x\to\pm\infty$). Therefore by Liouville  theorem, $G(z)$ has to be constant as $|z|\to\infty$, and so is $G(x)$ as $x\to\pm\infty$. Hence, the Lebesgue integral $$ \int_{-\infty}^{\infty} |e^{\lambda x^\alpha}|dx= \int_{-\infty}^{\infty} e^{\lambda_r x^\alpha} |e^{\lambda_i x^\alpha}|dx=\int_{-\infty}^{\infty} e^{\lambda_r x^\alpha}dx<\infty$$ since $G(x)$ is constant as $x\to\pm\infty$. For $\lambda_r<0$ and $\alpha$ odd, the limit $\lim_{x\to-\infty} e^{\lambda_r x^\alpha}$ diverges and so does the integral $\int_{-\infty}^{\infty} e^{\lambda_r x^\alpha}dx$. Therefore, the Lebesgue integral $ \int_{-\infty}^{\infty} |e^{\lambda x^\alpha}|dx$ has to diverge too. On the other hand, for $\lambda_r>0$, the limit $\lim_{x\to+\infty} e^{\lambda_r x^\alpha}$ diverges, and so does the integral  $\int_{-\infty}^{\infty} e^{\lambda_r x^\alpha}dx$ regardless of the value of $\alpha$. Therefore, the Lebesgue integral $ \int_{-\infty}^{\infty} |e^{\lambda x^\alpha}|dx$ has to diverge too.

\item At $x = 0,$ $g(0) = 1$. And so, around $x = 0$, the antiderivative $G(x) \sim x$
because $G^\prime(0) = g(0) = 1$.   And so  $(0, G(0)) = (0, 0)$. Moreover, $G^{\prime\prime}(x) =g^\prime(x) =\lambda\alpha x^{\alpha-1} e^{\lambda x^\alpha}, \alpha\ge 2$, gives $G^{\prime\prime}(0) = 0$. Hence, by the second derivative test, if $\lambda$ is real ($\lambda=\lambda_r$),  the point $(0, G(0)) = (0, 0)$ is an inflection
point of the curve $Y = G(x),$ $x \in\mathbb{ R}$.

\item For $\lambda=\lambda_r$ ($\lambda\in\mathbb{R}$), both $g(x)$ and $G(x)$ are analytic on $\mathbb{R}$. Using this fact and the fact that for even $\alpha$ and $\lambda_r<0$, $ \int_{-\infty}^{\infty} |e^{\lambda x^\alpha}|dx<\infty$  implies that for even $\alpha$ and $\lambda_r<0$, $G(x)$ has to be constant as $x\to\pm\infty$. In addition, the fact that $G^{\prime\prime}(x) < 0$ if $x < 0$ and $G^{\prime\prime}(x) > 0$ if $x > 0$ implies
that, $G(x)$ is concave upward on the interval $(−\infty, 0)$ while is concave downward on the interval $(0, +\infty)$. Moreover, the fact that $g(x) = G^\prime(x)$ is symmetric about the $y$-axis (even) implies that $G(x)$ has to be antisymmetric about the $y$-axis (odd). Hence there exists a real positive constant $\theta>0$ such that
 limits $\lim_{x\to-\infty}G(x){=}-\theta$ and $\lim_{x\to+\infty}G(x){=}\theta$.\hfill$\square$
\end{enumerate}
%\hfill$\square$\\[1ex]

\begin{figure}[t]
\centerline{\includegraphics[width=0.45\textwidth]{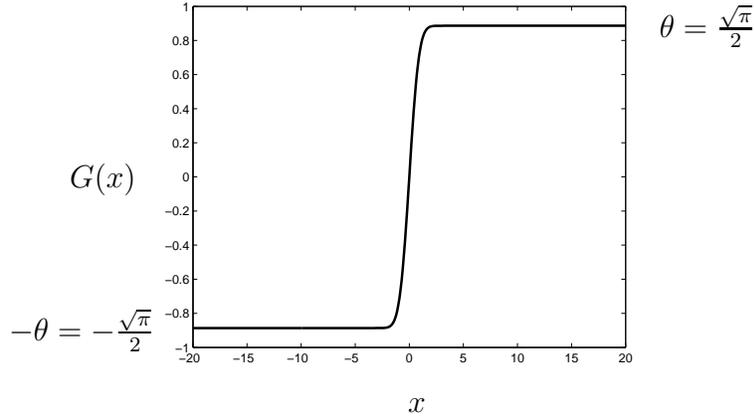}}
\begin{picture}(0,0)
\put (110,93){$G(x)$}
\put (238,8){$x$}
\put (87,35){$-\theta=-\frac{\sqrt{\pi}}{2}$}
\put (333,150){$\theta=\frac{\sqrt{\pi}}{2}$}
\end{picture}
\caption{$G(x)$ is the antiderivative of $e^{-x^2}$ given by (\ref{eq14}).}
\label{fig1}
\end{figure}

%\begin{exaen}
\emph{Example~3.}
If $\lambda=-1$ and $\alpha=2$, then
\begin{equation}
\int e^{-x^2}dx=x\ _1F_1\left(\frac{1}{2}; \frac{3}{2};  -x^2\right)+C.
\label{eq14}
\end{equation}
According to (\ref{eq14}), the antiderivative of $g(x) = e^{-x^2}$ is $G(x)=x\ _1F_1\left(\frac{1}{2}; \frac{3}{2};  -x^2\right)$. Its graph as a function of $x$, sketched using MATLAB, is shown in Figure \ref{fig1}. It is in agreement with Lemma \ref{lemen:1}. It is actually seen in Figure \ref{fig1} that $(0, 0)$ is an inflection point and that $G(x)$ reaches some constants as $x\to\pm\infty$ as predicted by Lemma \ref{lemen:1}.
\label{exaen:3}
%\end{exaen}

\medskip
In the following lemma, we obtain the values of $G(x)$, the antiderivative of the function $g(x) =e^{\lambda x^\alpha}$, as $x\to\pm\infty$ using the asymptotic expansion of the confluent hypergeometric function $_1F_1$.

\begin{lemen} Consider $G(x)$ in Proposition \ref{prpen:1}.
\begin{enumerate}
\item Then for $|x|\gg1$,
\begin{equation}
G(x)= x _1F_1\left(\frac{1}{\alpha}; \frac{1}{\alpha}+1; \lambda x^\alpha\right)\sim\left\{\begin{array}{l}
\Gamma\left( \frac{1}{\alpha}+1\right)\displaystyle\frac{e^{i\frac{\pi}{\alpha}}}{\lambda^{\frac{1}{\alpha}}}\frac{x}{|x|}+\frac{e^{\lambda x^\alpha}}{\alpha\lambda x^{\alpha-1}}, \mbox{if $\alpha$ is even}, \\[1.5ex]
 \Gamma\left( \frac{1}{\alpha}+1\right)\displaystyle\frac{e^{i\frac{\pi}{\alpha}}}{\lambda^{\frac{1}{\alpha}}}+\frac{e^{\lambda x^\alpha}}{\alpha\lambda x^{\alpha-1}}, \mbox{if $\alpha$ is odd}.\\
 \end{array}\right.
\label{eq14-1}
\end{equation}
\item Let $\alpha\ge2$ and be even, and let $\lambda=-\beta^2$, where $\beta$ is a real number, preferably positive. Then
\begin{equation}
G(-\infty)=\lim_{x\rightarrow-\infty}G(x)=\lim_{x\rightarrow-\infty} x _1F_1\left(\frac{1}{\alpha}; \frac{1}{\alpha}+1; -\beta^2 x^\alpha\right)=-\frac{1}{\beta^{\frac{2}{\alpha}}}\Gamma\left( \frac{1}{\alpha}+1\right)
\label{eq15}
\end{equation}
and
\begin{equation}
G(+\infty)=\lim_{x\rightarrow+\infty}G(x)=\lim_{x\rightarrow+\infty} x _1F_1\left(\frac{1}{\alpha}; \frac{1}{\alpha}+1; -\beta^2 x^\alpha\right)=\frac{1}{\beta^{\frac{2}{\alpha}}}\Gamma\left( \frac{1}{\alpha}+1\right).
\label{eq16}
\end{equation}
\item And by the FTC,
\begin{multline}
\int\limits_{-\infty}^{\infty} e^{-\beta^2x^\alpha} dx=G(+\infty)-G(-\infty)\\=\frac{1}{\beta^{\frac{2}{\alpha}}}\Gamma\left( \frac{1}{\alpha}+1\right)-\left(
-\frac{1}{\beta^{\frac{2}{\alpha}}}\Gamma\left( \frac{1}{\alpha}+1\right)\right)=\frac{2}{\beta^{\frac{2}{\alpha}}}\Gamma\left( \frac{1}{\alpha}+1\right).
\label{eq17}
\end{multline}
\end{enumerate}
\label{lemen:2}
\end{lemen}

\proofen
\begin{enumerate}
\item To prove (\ref{eq14-1}), we use the asymptotic series for the confluent hypergeometric function that is valid for $|z|\gg 1$ (\cite{1}, formula 13.5.1),
\begin{multline}
\frac{{}_1F_1\left(a;b;z\right)}{\Gamma(b)}=\frac{e^{\pm i\pi a}z^{-a}}{\Gamma(b-a)}\left\{\sum\limits_{n=0}^{R-1}\frac{(a)_n (1+a-b)_n}{n!}(-z)^{-n}+O(|z|^{-R})\right\}\\
+\frac{e^{z}z^{a-b}}{\Gamma(a)}\left\{\sum\limits_{n=0}^{S-1}\frac{(b-a)_n (1-a)_n}{n!}(z)^{-n}+O(|z|^{-S})\right\},
\label{eq18}
\end{multline}
where $a$ and $b$ are constants, and
the upper sign being taken if $-{\pi}/{2}<\text{arg}(z)<{3\pi}/{2}$ and the lower sign if $-{3\pi}/{2}<\text{arg}(z)\le -{\pi}/{2}$.
We set $z=\lambda x^\alpha, a=\frac{1}{\alpha}$ and $b=\frac{1}{\alpha}+1$, and obtain
\begin{multline}
\frac{{}_1F_1\left(\frac{1}{\alpha};\frac{1}{\alpha}+1;\lambda x^\alpha\right)}{\Gamma\left(\frac{1}{\alpha}+1\right)}=\frac{e^{i\frac{\pi}{\alpha}}}{(\lambda x^\alpha)^{\frac{1}{\alpha}}}\left\{\sum\limits_{n=0}^{R-1}\frac{\left(\frac{1}{\alpha}\right)_n }{n!}(\lambda x^\alpha)^{-n}+O\left\{\lambda x^\alpha\right)^{-R}\right\} \\
+\frac{e^{\lambda x^\alpha}(\lambda x^\alpha)^{-1}}{\Gamma\left(\frac{1}{\alpha}\right)}\left\{\sum\limits_{n=0}^{S-1} \left(1-\frac{1}{\alpha}\right)_n(\lambda x^\alpha)^{-n}+O\left(\lambda x^\alpha\right)^{-S}\right\}.
\label{eq19}
\end{multline}
Then, for $|x|\gg1$,
\begin{equation}
\frac{e^{i\frac{\pi}{\alpha}}}{(\lambda x^\alpha)^{\frac{1}{\alpha}}}\left\{\sum\limits_{n=0}^{R-1}\frac{\left(\frac{1}{\alpha}\right)_n }{n!}(\lambda x^\alpha)^{-n}+O\left\{\lambda x^\alpha\right)^{-R}\right\}\sim\left\{\begin{array}{l}
\displaystyle\frac{e^{i\frac{\pi}{\alpha}}}{\lambda^{\frac{1}{\alpha}}}\frac{1}{|x|}, \mbox{if $\alpha$ is even}, \\[1.5ex]
\displaystyle\frac{e^{i\frac{\pi}{\alpha}}}{\lambda^{\frac{1}{\alpha}}}\frac{1}{x}, \mbox{if $\alpha$ is odd},
 \end{array}\right.
\label{eq20}
\end{equation}
while
\begin{equation}
\frac{e^{\lambda x^\alpha}(\lambda x^\alpha)^{-1}}{\Gamma\left(\frac{1}{\alpha}\right)}\left\{\sum\limits_{n=0}^{S-1} \left(1-\frac{1}{\alpha}\right)_n(\lambda x^\alpha)^{-n}+O\left(\lambda x^\alpha\right)^{-S}\right\}\sim\frac{e^{\lambda x^\alpha}}{\Gamma\left(\frac{1}{\alpha}\right)\lambda x^{\alpha}}.
\label{eq21}
\end{equation}
And so, for $|x|\gg1$,
\begin{equation}
\frac{{}_1F_1\left(\frac{1}{\alpha};\frac{1}{\alpha}+1;\lambda x^\alpha\right)}{\Gamma\left( \frac{1}{\alpha}+1\right)}\sim\left\{\begin{array}{l}
\displaystyle\frac{e^{i\frac{\pi}{\alpha}}}{\lambda^{\frac{1}{\alpha}}}\frac{1}{|x|}+\frac{e^{\lambda x^\alpha}}{\Gamma\left(\frac{1}{\alpha}\right)\lambda x^{\alpha}}, \mbox{if $\alpha$ is even}, \\[1.5ex]
\displaystyle\frac{e^{i\frac{\pi}{\alpha}}}{\lambda^{\frac{1}{\alpha}}}\frac{1}{x}+\frac{e^{\lambda x^\alpha}}{\Gamma\left(\frac{1}{\alpha}\right)\lambda x^{\alpha}}, \mbox{if $\alpha$ is odd}.\\
 \end{array}\right.
\label{eq22}
\end{equation}
Hence,
\begin{equation}
G(x)= x _1F_1\left(\frac{1}{\alpha}; \frac{1}{\alpha}+1; \lambda x^\alpha\right)\sim\left\{\begin{array}{l}
\Gamma\left( \frac{1}{\alpha}+1\right)\displaystyle\frac{e^{i\frac{\pi}{\alpha}}}{\lambda^{\frac{1}{\alpha}}}\frac{x}{|x|}+\frac{e^{\lambda x^\alpha}}{\alpha\lambda x^{\alpha-1}}, \mbox{if $\alpha$ is even}, \\[1.5ex]
 \Gamma\left( \frac{1}{\alpha}+1\right)\displaystyle\frac{e^{i\frac{\pi}{\alpha}}}{\lambda^{\frac{1}{\alpha}}}+\frac{e^{\lambda x^\alpha}}{\alpha\lambda x^{\alpha-1}}, \mbox{if $\alpha$ is odd}.\\
 \end{array}\right.
\label{eq23}
\end{equation}

\item Setting $\lambda=-\beta^2$, where $\beta$ is real and positive and using (\ref{eq14-1}), then for $\alpha$
even,
\begin{equation}
G(x)= x _1F_1\left(\frac{1}{\alpha}; \frac{1}{\alpha}+1; -\beta^2 x^\alpha\right)\sim
\frac{1}{\beta^{\frac{2}{\alpha}}}\Gamma\left( \frac{1}{\alpha}+1\right)\frac{x}{|x|}-\frac{e^{-\beta^2 x^\alpha}}{\alpha\beta^2 x^{\alpha-1}}.
\label{eq24}
\end{equation}
Therefore,
\begin{equation}
G(-\infty)=\lim_{x\rightarrow-\infty}G(x)=\lim_{x\rightarrow-\infty} x _1F_1\left(\frac{1}{\alpha}; \frac{1}{\alpha}+1; -\beta^2 x^\alpha\right)=-\frac{1}{\beta^{\frac{2}{\alpha}}}\Gamma\left( \frac{1}{\alpha}+1\right)
\label{eq25}
\end{equation}
and
\begin{equation}
G(+\infty)=\lim_{x\rightarrow+\infty}G(x)=\lim_{x\rightarrow+\infty} x _1F_1\left(\frac{1}{\alpha}; \frac{1}{\alpha}+1; -\beta^2 x^\alpha\right)=\frac{1}{\beta^{\frac{2}{\alpha}}}\Gamma\left( \frac{1}{\alpha}+1\right).
\label{eq26}
\end{equation}

\item By the Fundamental Theorem of Calculus, we have
\begin{align}\begin{aligned}
\int\limits_{-\infty}^{+\infty} e^{-\beta^2x^\alpha} dx &=\lim_{y\rightarrow-\infty}\int\limits_{y}^{0} e^{-\beta^2x^\alpha} dx
+\lim_{y\rightarrow+\infty}\int\limits_{0}^{y} e^{-\beta^2x^\alpha}dx
\\ &=\lim_{y\rightarrow+\infty}y\ _1F_1\left(\frac{1}{\alpha};\frac{1}{\alpha}+1;-\beta^2y^\alpha\right)
-\lim_{y\rightarrow-\infty}y\ _1F_1\left(\frac{1}{\alpha};\frac{1}{\alpha}+1;-\beta^2y^\alpha\right)
\\ &
=G(+\infty)-G(-\infty)
\\ &=\frac{1}{\beta^{\frac{2}{\alpha}}}\Gamma\left( \frac{1}{\alpha}+1\right)-\left(
-\frac{1}{\beta^{\frac{2}{\alpha}}}\Gamma\left( \frac{1}{\alpha}+1\right)\right)=\frac{2}{\beta^{\frac{2}{\alpha}}}\Gamma\left( \frac{1}{\alpha}+1\right).
\label{eq26}
\end{aligned}\end{align}
\end{enumerate}

We now verify whether (\ref{eq26}) is correct or not by double integration. We first
observe that (\ref{eq26}) is valid for all even $\alpha\ge2$. And so, if (\ref{eq26}) is verified for
$\alpha=2$, we are done since (\ref{eq26}) is valid for all even $\alpha\ge2$. For $\alpha=2$,  we have
\begin{align}\begin{aligned}
\int\limits_{-\infty}^{+\infty} e^{-\beta^2x^2} dx &=\lim_{y\rightarrow-\infty}\int\limits_{y}^{0} e^{-\beta^2x^2} dx
+\lim_{y\rightarrow+\infty}\int\limits_{0}^{y} e^{-\beta^2x^2}dx
\\ &=\lim_{y\rightarrow+\infty}y\ _1F_1\left(\frac{1}{2};\frac{3}{2};-\beta^2y^2\right)
-\lim_{y\rightarrow-\infty}y\ _1F_1\left(\frac{1}{2};\frac{3}{2};-\beta^2y^2\right)
\\ &=G(+\infty)-G(-\infty)=\frac{2}{\beta}\Gamma\left(\frac{3}{2}\right)=\frac{2}{\beta}\frac{\sqrt{\pi}}{2}=\frac{\sqrt{\pi}}{\beta}.
\label{eq27}
\end{aligned}\end{align}

On the other hand,
\begin{align}
\left(\int\limits_{-\infty}^{\infty} e^{-\beta^2 x^2} dx\right)^2 &=\left(\int\limits_{-\infty}^{\infty} e^{-\beta^2 x^2} dx\right)\left(\int\limits_{-\infty}^{\infty} e^{-\beta^2 y^2} dy\right)
\\ &=\int\limits_{-\infty}^{\infty}\int\limits_{-\infty}^{\infty}e^{-\beta^2(x^2+y^2)}dydx.
\label{eq28}
\end{align}
In polar coordinate,
\begin{equation}
\int\limits_{-\infty}^{\infty}\int\limits_{-\infty}^{\infty}e^{-\beta^2(x^2+y^2)}dydx=\int\limits_{0}^{2\pi}\int\limits_{0}^{\infty}e^{-\beta^2 r^2}rdrd\theta
=\frac{1}{2\beta^2}\int\limits_{0}^{2\pi}d\theta=\frac{\pi}{\beta^2}.
\label{eq29}
\end{equation}
This gives
\begin{equation}
\int\limits_{-\infty}^{\infty} e^{-\beta^2x^2} dx=\sqrt{\int\limits_{-\infty}^{\infty}\int\limits_{-\infty}^{\infty}e^{-(x^2+y^2)}dydx}=\frac{\sqrt{\pi}}{\beta}
\label{eq30}
\end{equation}
as before.
\hfill$\square$\\[1ex]

%\begin{exaen}
\emph{Example~4.}
Setting $\lambda=-\beta^2=-1$, $\beta=1$ and $\alpha=2$ in Lemma \ref{lemen:2} gives
\begin{equation}
G(-\infty)=\lim_{x\rightarrow-\infty}G(x)=\lim_{x\rightarrow-\infty}x\ _1F_1\left(\frac{1}{2};\frac{3}{2};-x^2\right)
=-\frac{\sqrt{\pi}}{2}
\label{eq31}
\end{equation}
and
\begin{equation}
G(+\infty)=\lim_{x\rightarrow+\infty}G(x)=\lim_{x\rightarrow+\infty}x\ _1F_1\left(\frac{1}{2};\frac{3}{2};-x^2\right)
=\frac{\sqrt{\pi}}{2}.
\label{eq32}
\end{equation}
This implies $\theta={\sqrt{\pi}}/{2}$ in Lemma \ref{lemen:1}. And this is exactly the value of $G(x)$ as $x\to\infty$ in Figure \ref{fig1}. We also have $\lim_{x\rightarrow-\infty}G(x)=-\theta=-{\sqrt{\pi}}/{2}$ as in Figure \ref{fig1}. Using the FTC, we readily obtain
\begin{equation}
\int\limits_{-\infty}^{0} e^{-x^2} dx =G(0)-G(-\infty)=0-\left(-\frac{\sqrt{\pi}}{2}\right)=\frac{\sqrt{\pi}}{2},
\label{eq33}
\end{equation}
\begin{equation}
\int\limits_{0}^{+\infty} e^{-x^2} dx =G(+\infty)-G(0)=\frac{\sqrt{\pi}}{2}-0=\frac{\sqrt{\pi}}{2}
\label{eq34}
\end{equation}
and
\begin{equation}
\int\limits_{-\infty}^{+\infty} e^{-x^2} dx =G(+\infty)-G(-\infty)=\frac{\sqrt{\pi}}{2}-\left(-\frac{\sqrt{\pi}}{2}\right)=\sqrt{\pi}.
\label{eq35}
\end{equation}

\label{exaen:4}
%\end{exaen}
\medskip

%\begin{exaen}
\emph{Example~5.}
In this example, the integral

\begin{equation}
\int\limits_{-\infty}^{x} e^{t^{2n+1}} dt,\quad x<\infty,
\label{eq36}
\end{equation}
where $n$ is a positive integer, is evaluated using Proposition \ref{prpen:1} and the asymptotic expression  (\ref{eq14-1}). Setting $\lambda=1$ and $\alpha = 2n+ 1$ in Proposition \ref{prpen:1} , and using
 (\ref{eq14-1}) gives

\begin{align}\begin{aligned}
\int\limits_{-\infty}^{x} e^{t^{2n+1}} dt &=\lim_{y\rightarrow-\infty}\int\limits_{y}^{x} e^{t^{2n+1}} dt
\\ &=x\ _1F_1\left(\frac{1}{2n+1};\frac{2n+2}{2n+1};x^{2n+1}\right)
-\lim_{y\rightarrow-\infty}y\ _1F_1\left(\frac{1}{2n+1};\frac{2n+2}{2n+1};y^{2n+1}\right)
\\ &=x\ _1F_1\left(\frac{1}{2n+1};\frac{2n+2}{2n+1};x^{2n+1}\right)
-\Gamma\left(\frac{2n+2}{2n+1}\right),\quad x<\infty.
\label{eq37}
\end{aligned}\end{align}
One can also obtain
\begin{align}\begin{aligned}
\int\limits_{x}^{+\infty} e^{-t^{2n+1}} dt &=\lim_{y\rightarrow+\infty}\int\limits_{x}^{y} e^{-t^{2n+1}} dt
\\ &=\lim_{y\rightarrow-\infty}y\ _1F_1\left(\frac{1}{2n+1};\frac{2n+2}{2n+1};-y^{2n+1}\right)
-x\ _1F_1\left(\frac{1}{2n+1};\frac{2n+2}{2n+1};-x^{2n+1}\right)
\\ &=\Gamma\left(\frac{2n+2}{2n+1}\right)-x\ _1F_1\left(\frac{1}{2n+1};\frac{2n+2}{2n+1};-x^{2n+1}\right)
,\quad x>-\infty.
\label{eq38}
\end{aligned}\end{align}
\label{exaen:5}
%\end{exaen}

\begin{teoen} For any $A$ and $B$, the FTC gives
\begin{equation}
\int\limits_A^B e^{\lambda x^\alpha}dx = G(B)-G(A),
\label{eq39}
\end{equation}
where $G$ is the antiderivative of the function $g(x) = e^{\lambda x^\alpha}$
and is given in Proposition \ref{prpen:1}. And $\lambda$ is any complex or real constant, and $\alpha\ge2$.
\label{teoen:1}
\end{teoen}

\proofen
$G(x) = x\ _1F_1\left(\frac{1}{\alpha}; \frac{1}{\alpha}+1; \lambda x^\alpha\right)$, where $\lambda$ is any constant, is the antiderivative of  $g(x) = e^{\lambda x^\alpha}, \alpha\ge2$ by Proposition \ref{prpen:1}, Lemma \ref{lemen:1} and Lemma \ref{lemen:2}. And since the FTC works for $A =-\infty$ and $B = 0$ in (\ref{eq33}), $A = 0$ and $B = +\infty$ in (\ref{eq34}) and $A =-\infty$ and $B = +\infty$ in (\ref{eq35}) by Lemma 2 if $\lambda = −1$ and
$\alpha = 2$, and for all $\lambda < 0$ and all even $\alpha \ge 2$, then it has to work for other
values of $A, B \in \mathbb{R}$ and for any $\lambda\in\mathbb{C}$ and $\alpha \ge 2$. This completes the proof.
\hfill$\square$\\[1ex]

\emph{Example~6.} In this example, we apply Theorem \ref{teoen:1} to the Central Limit Theorem
in Probability theory \cite{2}. The normal zero-one distribution of a random
variable X is the measure $\mu(dx) = g_X(x)dx$, where $dx$ is the Lebesgue measure
and the function $g_X(x)$ is the probability density function (p.d.f) of the
normal zero-one distribution \emph{\cite{2}}, and is
\begin{equation}
g_{X}(x) =\frac{1}{\sqrt{2\pi}}e^{-\frac{x^2}{2}}, -\infty<x<+\infty.
\label{eq40}
\end{equation}
A comparison with the function g(x) in Proposition \ref{prpen:1} and Lemma \ref{lemen:1} gives $\lambda = −\beta^2=-1/2$ and $\alpha = 2$. By Theorem 1, the cumulative probability, $P(X < z)$, is then given by

\begin{equation}
P(X<z)=\mu\{(-\infty,z)\}=\int\limits_{-\infty}^z g_X(x) dx=\frac{1}{\sqrt{2\pi}}\int\limits_{-\infty}^z e^{-\frac{x^2}{2}}dx=\frac{1}{2}+\frac{z}{\sqrt{2\pi}}\ _1F_1\left(\frac{1}{2};\frac{3}{2};-\frac{z^2}{2}\right).
\label{eq41}
\end{equation}
\label{exaen:6}
For example, we can also use Theorem \ref{teoen:1} to obtain $P(-2 < X < 2) =
\mu{(-2, 2)} = 0.4772 - (-0.4772) = 0.9544$, $P(-1 < X < 2) = \mu{(-1, 2)} =0.4772 - (-0.3413) = 0.8185$ and so on.
%\end{exaen}

\medskip

\emph{Example~7.} Using integration by parts and applying Theorem \ref{teoen:1}, the Maxwell-Bortsman distribution is written in terms of the confluent hypergeometric $_1F_1$ as
\begin{equation}
F(v)=\theta\int\limits_{0}^v  x^2e^{-\gamma x^2}dx=-\frac{\theta v}{2\gamma}e^{-\gamma v^2}+\frac{\theta v}{2\gamma}\ _1F_1\left(\frac{1}{2}; \frac{3}{2}; -\gamma v^2\right)=\frac{\theta v}{2\gamma}\left[_1F_1\left(\frac{1}{2}; \frac{3}{2}; -\gamma v^2\right)-e^{-\gamma v^2}\right].
\label{eq42}
\end{equation}
\label{exaen:7}

\section{Other related non-elementary integrals} \label{sec:3}
\setcounter{equation}{0}

\begin{prpen}The function $G(x) = x\ _1F_2\left(\frac{1}{2\alpha}; \frac{1}{2},\frac{1}{2\alpha}+1; \frac{\lambda^2 x^{2\alpha}}{4}\right)$, where $_1F_2$ is a hypergeometric function \emph{\cite{1}}, $\lambda$ is an arbitrarily constant and $\alpha\ge 2$,
is the antiderivative of the function $g(x) = \cosh{(\lambda x^\alpha)}$. Thus,
\begin{equation}
\int \cosh{(\lambda x^\alpha)}dx=x\ _1F_2\left(\frac{1}{2\alpha}; \frac{1}{2},\frac{1}{2\alpha}+1; \frac{\lambda^2 x^{2\alpha}}{4}\right)+C.
\label{eq43}
\end{equation}
\label{prpen:2}
\end{prpen}

\proofen   We proceed as before. We expand $g(x) = \cosh{(\lambda x^\alpha)}$ as a Taylor series and integrate the series
term by term, use the Pochhammer’s notation \cite{1} for the gamma
function, $\Gamma(a + n) = \Gamma(a)(a)_n$, where $(a)_n = a(a + 1) \cdots (a + n - 1)$,
and the property of the gamma function $\Gamma(a + 1) = a\Gamma(a)$ \cite{1}. We also use the Gamma duplication formula \cite{1}. We then obtain
\begin{align}\begin{aligned}
\int g(x) dx &=\int \cosh{(\lambda x^\alpha)}dx=\sum\limits_{n=0}^{\infty}\frac{\lambda^{2n}}{(2n)!}\int x^{2\alpha n}dx  \\&
=\sum\limits_{n=0}^{\infty}\frac{\lambda^{2n}}{(2n)!}\frac{x^{2\alpha n+1}}{2\alpha n+1}+C
 \\&=\frac{x}{2\alpha}\sum\limits_{n=0}^{\infty}\frac{(\lambda^2 x^{2\alpha})^n}{(2n)!\left(n+\frac{1}{2\alpha}\right)}+C \\&=\frac{x}{2\alpha}\sum\limits_{n=0}^{\infty}\frac{\Gamma\left(n+\frac{1}{2\alpha}\right)}{\Gamma(2n+1)\Gamma\left(n+\frac{1}{2\alpha}+1\right)}{(\lambda^2 x^{2\alpha})^n}+C
 \\&
={x}\sum\limits_{n=0}^{\infty}\frac{\left(\frac{1}{2\alpha}\right)_n}{\left(\frac{1}{2}\right)_n\left(\frac{1}{2\alpha}+1\right)_n}\frac{(\lambda^2 x^{2\alpha})^n}{n!}+C
\\&=x\ _1F_2\left(\frac{1}{2\alpha}; \frac{1}{2},\frac{1}{2\alpha}+1; \frac{\lambda^2 x^{2\alpha}}{4}\right)+C=G(x) + C.\quad \square
\label{eq44}
\end{aligned}\end{align}

%  \hfill$\square$\\[1ex]%--- P r o o f.

\begin{prpen}The function $$G(x) =\frac{\lambda x^{\alpha+1}}{\alpha+1}\ _1F_2\left(\frac{1}{2\alpha}+\frac{1}{2}; \frac{3}{2},\frac{1}{2\alpha}+ \frac{3}{2}; \frac{\lambda^2 x^{2\alpha}}{4}\right),$$ where $_1F_2$ is a hypergeometric function \emph{\cite{1}}, $\lambda$ is an arbitrarily constant and $\alpha\ge 2$,
is the antiderivative of the function $g(x) = \sinh{(\lambda x^\alpha)}$. Thus,
\begin{equation}
\int \sinh{(\lambda x^\alpha)}dx=\frac{\lambda x^{\alpha+1}}{\alpha+1}\ _1F_2\left(\frac{1}{2\alpha}+\frac{1}{2}; \frac{3}{2},\frac{1}{2\alpha}+ \frac{3}{2}; \frac{\lambda^2 x^{2\alpha}}{4}\right)+C.
\label{eq45}
\end{equation}
\label{prpen:3}
\end{prpen}

\proofen   As above, we expand $g(x) = \sinh{(\lambda x^\alpha)}$ as a Taylor series and integrate the series
term by term, use the Pochhammer’s notation \cite{1} for the gamma
function, $\Gamma(a + n) = \Gamma(a)(a)_n$, where $(a)_n = a(a + 1) \cdots (a + n - 1)$,
and the property of the gamma function $\Gamma(a + 1) = a\Gamma(a)$ \cite{1}. We also use the Gamma duplication formula~\cite{1}. We then obtain
\begin{align}\begin{aligned}
\int g(x) dx &=\int \sinh{(\lambda x^\alpha)}dx=\sum\limits_{n=0}^{\infty}\frac{\lambda^{2n+1}}{(2n+1)!}\int x^{2\alpha n+\alpha}dx \\&
=\sum\limits_{n=0}^{\infty}\frac{\lambda^{2n+1}}{(2n+1)!}\frac{x^{2\alpha n+\alpha+1}}{2\alpha n+\alpha+1}+C
 \\&=\frac{\lambda x^{\alpha+1}}{2\alpha}\sum\limits_{n=0}^{\infty}\frac{(\lambda^2 x^{2\alpha})^n}{(2n+1)!\left(n+\frac{1}{2\alpha}+\frac{1}{2}\right)}+C \\&=\frac{\lambda x^{\alpha+1}}{2\alpha}\sum\limits_{n=0}^{\infty}\frac{\Gamma\left(n+\frac{1}{2\alpha}+\frac{1}{2}\right)}{\Gamma(2n+2)\Gamma\left(n+\frac{1}{2\alpha}
 +\frac{3}{2}\right)}{(\lambda^2 x^{2\alpha})^n}+C
r \\&
=\frac{\lambda x^{\alpha+1}}{\alpha+1}\sum\limits_{n=0}^{\infty}\frac{\left(\frac{1}{2\alpha}+
\frac{1}{2}\right)_n}{\left(\frac{3}{2}\right)_n\left(\frac{1}{2\alpha}+\frac{3}{2}\right)_n}\frac{(\lambda^2 x^{2\alpha})^n}{n!}+C
 \\&=\frac{\lambda x^{\alpha+1}}{\alpha+1}\ _1F_2\left(\frac{1}{2\alpha}+\frac{1}{2}; \frac{3}{2},\frac{1}{2\alpha}+ \frac{3}{2}; \frac{\lambda^2 x^{2\alpha}}{4}\right)+C=G(x) + C.\quad \square
\label{eq46}
\end{aligned}\end{align}
%  \hfill$\square$\\[1ex]%--- P r o o f.
We also can show as above that
\begin{equation}
\int \cos{(\lambda x^\alpha)}dx=x\ _1F_2\left(\frac{1}{2\alpha}; \frac{1}{2},\frac{1}{2\alpha}+1; -\frac{\lambda^2 x^{2\alpha}}{4}\right)+C
\label{eq47}
\end{equation}
and
\begin{equation}
\int \sin{(\lambda x^\alpha)}dx=\frac{\lambda x^{\alpha+1}}{\alpha+1}\ _1F_2\left(\frac{1}{2\alpha}+\frac{1}{2}; \frac{3}{2},\frac{1}{2\alpha}+ \frac{3}{2}; -\frac{\lambda^2 x^{2\alpha}}{4}\right)+C.
\label{eq48}
\end{equation}

\medskip
\begin{teoen}
For any constants $\alpha$ and $\lambda$,
\begin{equation}
 _1F_2\left(\frac{1}{2\alpha}; \frac{1}{2},\frac{1}{2\alpha}+1; \frac{\lambda^2 x^{2\alpha}}{4}\right)\\=\frac{1}{2}\left[ _1F_1\left(\frac{1}{\alpha}; \frac{1}{\alpha}+1; \lambda x^\alpha\right)+\ _1F_1\left(\frac{1}{\alpha}; \frac{1}{\alpha}+1; -\lambda x^\alpha\right)\right]
\label{eq49}
\end{equation}
and
\begin{equation}
 _1F_2\left(\frac{1}{2\alpha}; \frac{1}{2},\frac{1}{2\alpha}+1; -\frac{\lambda^2 x^{2\alpha}}{4}\right)\\=\frac{1}{2}\left[ _1F_1\left(\frac{1}{\alpha}; \frac{1}{\alpha}+1; i\lambda x^\alpha\right)+\ _1F_1\left(\frac{1}{\alpha}; \frac{1}{\alpha}+1; -i\lambda x^\alpha\right)\right].
\label{eq50}
\end{equation}

\label{teoen:2}
\end{teoen}

\proofen  Using Proposition \ref{prpen:1}, we obtain
\begin{multline}
\int \cosh{(\lambda x^\alpha)}dx=\int \frac{e^{\lambda x^\alpha}+ e^{-\lambda x^\alpha}}{2}dx\\=\frac{x}{2}\left[ _1F_1\left(\frac{1}{\alpha}; \frac{1}{\alpha}+1; \lambda x^\alpha\right)+\ _1F_1\left(\frac{1}{\alpha}; \frac{1}{\alpha}+1; -\lambda x^\alpha\right)\right]+ C.
\label{eq51}
\end{multline}
Hence, comparing (\ref{eq43}) with (\ref{eq51}) gives (\ref{eq49}). Using Proposition \ref{prpen:1}, on the other hand, we obtain
\begin{multline}
\int \cos{(\lambda x^\alpha)}dx=\int \frac{e^{i\lambda x^\alpha}+ e^{-i\lambda x^\alpha}}{2}dx\\=\frac{x}{2}\left[ _1F_1\left(\frac{1}{\alpha}; \frac{1}{\alpha}+1; i\lambda x^\alpha\right)+\ _1F_1\left(\frac{1}{\alpha}; \frac{1}{\alpha}+1; -i\lambda x^\alpha\right)\right]+ C.
\label{eq52}
\end{multline}
Hence, comparing (\ref{eq47}) with (\ref{eq52}) gives (\ref{eq50}).
  \hfill$\square$\\[1ex]%--- P r o o f.

\begin{teoen}
For any constants $\alpha$ and $\lambda$,

\begin{multline}
 \frac{\lambda x^{\alpha}}{\alpha+1}\ _1F_2\left(\frac{1}{2\alpha}+\frac{1}{2}; \frac{3}{2},\frac{1}{2\alpha}+ \frac{3}{2}; -\frac{\lambda^2 x^{2\alpha}}{4}\right)\\=\frac{1}{2}\left[ _1F_1\left(\frac{1}{\alpha}; \frac{1}{\alpha}+1; \lambda x^\alpha\right)-\ _1F_1\left(\frac{1}{\alpha}; \frac{1}{\alpha}+1; -\lambda x^\alpha\right)\right]
\label{eq53}
\end{multline}
and
\begin{multline}
\frac{\lambda x^{\alpha}}{\alpha+1}\ _1F_2\left(\frac{1}{2\alpha}+\frac{1}{2}; \frac{3}{2},\frac{1}{2\alpha}+ \frac{3}{2}; -\frac{\lambda^2 x^{2\alpha}}{4}\right)\\=\frac{1}{2i}\left[ _1F_1\left(\frac{1}{\alpha}; \frac{1}{\alpha}+1; i\lambda x^\alpha\right)-\ _1F_1\left(\frac{1}{\alpha}; \frac{1}{\alpha}+1; -i\lambda x^\alpha\right)\right].
\label{eq54}
\end{multline}

\label{teoen:3}
\end{teoen}

\proofen  Using Proposition \ref{prpen:1}, we obtain
\begin{multline}
\int \sinh{(\lambda x^\alpha)}dx=\int \frac{e^{\lambda x^\alpha}+ e^{-\lambda x^\alpha}}{2}dx\\=\frac{x}{2}\left[ _1F_1\left(\frac{1}{\alpha}; \frac{1}{\alpha}+1; \lambda x^\alpha\right)-\ _1F_1\left(\frac{1}{\alpha}; \frac{1}{\alpha}+1; -\lambda x^\alpha\right)\right]+ C.
\label{eq55}
\end{multline}
Hence, comparing (\ref{eq45}) with (\ref{eq55}) gives (\ref{eq53}). Using Proposition \ref{prpen:1}, on the other hand, we obtain
\begin{multline}
\int \sin{(\lambda x^\alpha)}dx=\int \frac{e^{i\lambda x^\alpha}+ e^{-i\lambda x^\alpha}}{2i}dx\\=\frac{x}{2i}\left[ _1F_1\left(\frac{1}{\alpha}; \frac{1}{\alpha}+1; i\lambda x^\alpha\right)-\ _1F_1\left(\frac{1}{\alpha}; \frac{1}{\alpha}+1; -i\lambda x^\alpha\right)\right]+ C.
\label{eq56}
\end{multline}
Hence, comparing (\ref{eq48}) with (\ref{eq56}) gives (\ref{eq54}).
  \hfill$\square$\\[1ex]%--- P r o o f.

\section{Conclusion}\label{sec:4}
\setcounter{equation}{0}

The non-elementary integral $\int e^{\lambda x^\alpha}dx$, where $\lambda$ is an arbitrary constant and $\alpha\ge 2$,
was expressed in term of the confluent hypergeometric function $_1F_1$. And using the properties of the confluent hypergeometric function $_1F_1$, the asymptotic expression for $|x|\gg 1$ of this integral was derived too. As established in
Theorem \ref{teoen:1}, the definite integral (\ref{eq6}) can now be computed using the FTC. For example, one can evaluate the area under the Gaussian Bell curve using the FTC rather than using double integration and then polar coordinates. One can also choose to use Theorem 1 to compute the cumulative probability for the normal distribution or that for the Maxwell-Bortsman distribution as shown in examples \ref{exaen:6} and \ref{exaen:7}.

On one hand, the integrals $\int \cosh(\lambda x^\alpha) dx,$ $\int \sinh(\lambda x^\alpha) dx,$ $\int \cos(\lambda x^\alpha) dx$ and $\int \sin(\lambda x^\alpha) dx,$ ${\alpha \ge 2}$, were evaluated in terms of the confluent hypergeometric function $_1F_1$, while on another hand, they were expressed in terms of the hypergeometric $_1F_2$. This allowed to express the hypergeometric function$_1F_2$ in terms of the confluent hypergeometric function $_1F_1$ (Theorems \ref{teoen:2} and~\ref{teoen:3}).

\begin{Biblioen}

\bibitem{1}{\bf Abramowitz~M., Stegun~I.A.} Handbook of mathematical functions with formulas,
graphs and mathematical tables. National Bureau of Standards, 1964. 1046~p.
\bibitem{2}{\bf Billingsley~P.} Probability and measure. Wiley series in Probability and Mathematical
Statistics, 3rd~Edition, 1995. 608~p.
\bibitem{3}{\bf Krantz~S.G.} Handbook of complex variables. Boston: MA Birkh\"{a}user, 1999. 290~p. DOI:~10.1007/978-1-4612-1588-2
\bibitem{4}{\bf Marchisotto~E.A., Zakeri~G.-A.} An invitation to integration in finite terms // College
Math. J., 1994. Vol.~25, no.~4. P.~295--308. DOI:~10.2307/2687614
\bibitem{5} NIST Digital Library of Mathematical Functions. \url{http://dlmf.nist.gov/}
\bibitem{6}{\bf Rosenlicht~M.} Integration in finite terms // Amer. Math. Monthly, 1972. Vol~79, no.~9. P.~963--972.
DOI:~10.2307/2318066
\end{Biblioen}

 \end{document}